\documentclass[11pt]{article}
\usepackage{amsthm}
\usepackage{latexsym, amssymb, amsmath}
\usepackage{graphicx}
\usepackage[dvipsnames]{xcolor}

\begin{document}
\numberwithin{equation}{section}

%%%%%%%%%%%%%%%%%%%%%%%%%%%%%%%%%%%%%%%%%%%%%%%%%%%%%%%%%%%%%%%%
\newtheorem{THEOREM}{Theorem}
\newtheorem{PRO}{Proposition}
\newtheorem{XXXX}{\underline{Theorem}}
\newtheorem{CLAIM}{Claim}
\newtheorem{COR}{Corollary}
\newtheorem{LEMMA}{Lemma}
\newtheorem{REM}{Remark}
\newtheorem{EX}{Example}
\newenvironment{PROOF}{{\bf Proof}.}{{\ \vrule height7pt width4pt depth1pt} \par \vspace{2ex} }
\newcommand{\Bibitem}[1]{\bibitem{#1} \ifnum\thelabelflag=1 
  \marginpar{\vspace{0.6\baselineskip}\hspace{-1.08\textwidth}\fbox{\rm#1}}
  \fi}
\newcounter{labelflag} \setcounter{labelflag}{0}
\newcommand{\labelon}{\setcounter{labelflag}{1}}
\newcommand{\Label}[1]{\label{#1} \ifnum\thelabelflag=1 
  \ifmmode  \makebox[0in][l]{\qquad\fbox{\rm#1}}
  \else\marginpar{\vspace{0.7\baselineskip}\hspace{-1.15\textwidth}\fbox{\rm#1}}
  \fi \fi}
% \labelon

\newcommand{\LEFTLINE}{\ifhmode\newline\else\noindent\fi}
\newcommand{\RIGHTLINE}[1]{\LEFTLINE\rightline{#1}}
\newcommand{\CENTERLINE}[1]{\LEFTLINE\centerline{#1}}
\def\BOX #1 #2 {\framebox[#1in]{\parbox{#1in}{\vspace{#2in}}}}
\parskip=8pt plus 2pt
\def\AUTHOR#1{\author{#1} \maketitle}
\def\Title#1{\begin{center}  \Large\bf #1 \end{center}  \vskip 1ex }
\def\Author#1{\vspace*{-2ex}\begin{center} #1 \end{center}  
 \vskip 2ex \par}
\renewcommand{\theequation}{\arabic{section}.\arabic{equation}}
\def\bdk#1{\makebox[0pt][l]{#1}\hspace*{0.03ex}\makebox[0pt][l]{#1}\hspace*{0.03ex}\makebox[0pt][l]{#1}\hspace*{0.03ex}\makebox[0pt][l]{#1}\mbox{#1} }
\def\psbx#1 #2 {\mbox{\psfig{file=#1,height=#2}}}

%pm  /mo  monotone
%pm  /po  polynomial
%pm  /py  positivity
%pm  /PS  partial sum
%pm  /ni  non-increasing
%pm  /tg  trigonometric
%pm  /br  \fbox{ \begin{minipage}[t]{5.4in}
%pm  /er  \end{minipage} }
%pm  /HR  {\par\vspace*{\baselineskip}\hrule{\vspace*{\baselineskip}\par}
%pm  /HO  H\^opital
%pm  /LH  L'H\^opital
%pm  /nd  non-decreasing
 
\newcommand{\ia}{\,\,\nearrow}
\newcommand{\da}{\,\,\searrow}
\newcommand{\ic}{\nearrow}
\newcommand{\dc}{\searrow}
\newcommand{\me}{\mbox{e}}
\newcommand{\vp}[2]{\vphantom{\vrule height#1pt depth#2pt}}
\newcommand{\FG}[2]{{\includegraphics[height=#1mm]{#2.eps}}}
\renewcommand{\thefootnote}{\fnsymbol{footnote}}
 
%%%%%%%%%%%%%%%%%%%%%%%%%%%%%%%%%%%%%%%%%%%%%%%%%%%%%%%%%%%%%%%%
 
\renewcommand{\theequation}{\arabic{equation}}
\Title{On H\^opital-style rules for monotonicity and oscillation}

% .a Man Kam Kwong

\vspace{0.5cm}
\begin{center}
MAN KAM KWONG\footnote{The research of this author is supported by the Hong Kong Government GRF Grant PolyU 5003/12P and the Hong Kong Polytechnic University Grants G-UC22 and G-YBCQ.}
\end{center}

% \vspace{0.5cm}
\begin{center}
\emph{Department of Applied Mathematics\\ The Hong Kong Polytechnic University,\\ Hunghom, Hong Kong}\\
\tt{mankwong@polyu.edu.hk}
\end{center}

\par\vspace*{\baselineskip}\par

\newcommand{\Cr}{\color{red}}
\newcommand{\mb}{\mathbf}
\newcommand{\tx}[1]{\mbox{ #1 }}

\parskip=6pt

\begin{abstract}
We point out the connection of the so-called H\^opital-style rules
for monotonicity and oscillation
to some well-known properties of concave{/}convex functions.
From this standpoint, we are able to generalize
the rules under no differentiability requirements 
and greatly extend their usability.
The improved rules can handle situations in which the functions 
involved have non-zero
initial values and when the derived functions are not necessarily monotone.
This perspective is not new; it can be dated back to Hardy, Littlewood
and P\'olya. 
\end{abstract}

\vspace{0.9cm}

\vspace{0.2cm}

\newpage
\section{Introduction and historical remarks}

Since the 1990's, many authors have successfully applied the so-called
monotone L'H\^opital's\footnote{A well-known anecdote, recounted in some undergraduate
textbooks and the Wikipedia, claims that H\^opital might have
cheated his teacher Johann Bernoulli to earn the credit for this classical
rule.}
rules to establish many new and useful inequalities.

In this article we make use of the connection of these rules
to some well-known properties of concave{/}convex functions (via a change of variable)
to extend the usability of the rules. 
We show how
the rules can be formulated under no differentiability requirements,
and how to characterize all possible situations, which are 
normally not covered by the conventional form of the rules.
These include situations in which the functions involved have non-zero
initial values and when the derived functions are not necessarily monotone.
This perspective is not new; it can be dated back to Hardy, Littlewood
and P\'olya (abbreviated as HLP). 

The concepts of non-decreasing, increasing (we use this term in the sense of what
some authors prefer to call ``strictly increasing''), non-increasing, and decreasing functions are
defined as usual. The term ``monotone'' can refer to any
one of these senses. For convenience, we use the symbols 
$ \nearrow $ and $ \searrow $ to denote increasing and decreasing, respectively.
The rules have appeared in various formulations. Let us 
start with the most popularly known form. Let $ a<b\leq \infty  $.

\par\vspace*{-3mm}\par
\begin{center}
\fbox{
\begin{minipage}[t]{5in}
\em 
\par\vspace*{2mm}\par
Let $ f,g:[a,b)\rightarrow \mathbb R $ be two continuous real-valued functions,
such that $ f(a)=g(a)=0 $, and $ g(x)>0 $ for $ x>a $.
Assume that they are differentiable at each point in $ (a,b) $, and 
$ g'(x)>0 $ for $ x>a $.
\par\vspace*{2mm}\par
\CENTERLINE{If \,\, $\displaystyle \frac{f'(x)}{g'(x)} $ \,\, is monotone in $ [a,b) $, so is \,\, $\displaystyle \frac{f(x)}{g(x)} $ ``in the same sense''.}
\par\vspace*{2mm}\par
\end{minipage}
}
\end{center}

A dual form assumes $ b<\infty  $, $ f(b)=g(b)=0 $ and $ g'(x)<0 $.
The same conclusion holds.

Note that even though the rule allows 
$ b=\infty  $, and that $ f $ or $ g $ need not be defined at $ b $, as far as the proof 
is concerned we may
assume without loss of generality that $ b $ is finite, and $ f $ and $ g $ are defined and continuous up to $ b $,
because we can first study the functions in a smaller subinterval 
$ [\alpha ,\beta ]\subset(0,b) $ and then let $ \alpha \rightarrow a $ and $ \beta \rightarrow b $.
When $ f $ and $ g $ are assumed to be differentiable, we take $ f'(a) $ and 
$ g'(a) $ to mean the righthand derivative at $ x=a $, and $ f'(b) $, $ g'(b) $ to mean
the lefthand derivative at $ x=b $.
Moreover, we can take $ a=0 $ after a suitable translation. 
These simplifications will be assumed in the rest of the paper.

Some of the variations in other formulations are merely cosmetic. For
example, if $ f(0)\neq 0 $, then use $ f(x)-f(0) $ instead, or if $ f $ is not defined
at $ 0 $, use the righthand limit $ f(0+) $, if it exists. We will not dwell 
further on these types. 

Some other variations are attempts to 
weaken the differentiability requirement on $ f $ and~ $ g $. A different type
concerns stronger formulations in which strict monotonicity can 
be deduced from
non-strict hypothesis. We defer a discussion of such variations to Section~2.
In the majority of concrete practical applications, however,
the functions involved
are fairly smooth, often infinitely differentiable and the stronger
forms are seldom used.

Propositions 147 and 148 (page 106) in the famous classic by HLP
\cite{hlp} (First Edition 1934) read:

\par\vspace*{-3mm}\par
\begin{center}
\fbox{
\begin{minipage}[t]{5in}
\em 
\par\vspace*{2mm}\par
{\bf 147.} The function
\par\vspace*{-6mm}\par
$$  \sigma (x) = \frac{\vp{17}{13}\displaystyle \int_{0}^{x} (1+ \sec t)\,\log\,\sec t \,dt }{\vp{18}{16}\displaystyle \log\,\sec x \int_{0}^{x} (1+\sec t)\,dt}  $$
increases steadily from $ \frac{1}{3} $ to $ \frac{1}{2} $ as $ x $ increases from $ 0 $ to $ \frac{1}{2} \pi  $.

\par\vspace*{\baselineskip}\par
There is a general theorem which will be found useful in the proof of Theorem 147.
\par\vspace*{\baselineskip}\par

{\bf 148.} If $ f $, $ g $, and $ f'/g' $ are {\color{red}positive increasing} functions, then $ f/g $ either increases
for all $ x $ in question, {\Cr or decreases for all such $ x $, or decreases to 
   a minimum and then increases.} In particular, if $ (0)=g(0)=0,then $ $ f/g $
   increases for $ x>0 $.
\par\vspace*{2mm}\par
\end{minipage}
}
\end{center}

Other than the phrases in red, Proposition 148 is essentially the 
increasing part of the monotone
rule, and Proposition 147 is an application of the rule in the same 
spirit as in the more recent work. At first reading, Proposition 148
appears to be weaker than the modern rule because it requires the 
additional conditions that $ f $ is + and $ \nearrow $, and $ f'/g' $ is +. 
Let us take a closer look at HLP's short and
elegant proof which is reproduced below.

\begin{PROOF}
[\,Hardy, Littlewood, P\'olya\,] To prove this, observe that
$$  \frac{d}{dx} \left( \frac{f}{g} \right)  = \left( \frac{f'}{g'} - \frac{f}{g} \right) \,\frac{g'}{g}  $$
and consider the possible intersections of the curves $ y=f/g $, $ y=f'/g' $.
At one of these intersections the first curve has a horizontal and the
second a rising tangent, and therefore there can be at most one intersection.

If we take $ g $ as the independent variable, write $ f(x)=\phi (g) $, and suppose,
as in the last clause of the theorem, that
$$  f(0)=g(0)=0,  $$
or $ \phi (0)=0 $, then the theorem takes the form: {\em if $ \phi (0)=0 $ and $ \phi '(g) $
increases for $ g>0 $, then $ \phi /g $ increases for $ g>0 $.} This is a slight
generalization of part of Theorem 127.
\end{PROOF}

In the proof, the + $ \ic $ property
of $ g $ is needed to guarantee that the denominators of the fractions
$ f/g $ and $ f'/g' $ will not become 0. The + $ \ic $ property of
$ f $ and the positivity of $ f'/g' $, however, is not needed anywhere.

Once the extra conditions are disposed of, we see that the $ \dc $ version of rule 
also holds, by considering $ -f(x) $ instead of $ f(x) $.

HLP did not attribute the result to anyone. It could mean that it is
one of their own, or it was widely known. We did not attempt to track
it down further in earlier literature.

In one of the Bourbaki books, \cite{bo} (1958), Exercise  10 of Chapter 1, \S2
(page 38) is exactly HLP's Proposition 148, minus those extra
conditions. It seems that in the intervening years,
someone must have figured out that those are superfluous. 
Neither a solution nor any attribution is given.
Mitrinovic \cite{mi} quoted Bourbaki's Exercise as \S3.9.49.

HLP's Proposition, on the other hand, is broader than the modern rule.
It says something about the general situation when neither $ f(0)=0 $
nor $ g(0)=0 $ is assumed. In Section~3, we will fully characterize all such
situations.

A classroom note by Mott \cite{mo} (1963) in the Monthly
presented a much weaker version of the monotone rule, in the
integral formulation (see Section~2). In the (more or less)
equivalent differential formulation,
the hypotheses require $ f'\ic $ and $ g'\dc $  (which implies
$ f'/g' $ is $ \ic $). Two follow-up papers by Redheffer \cite{re} (1964) and 
Boas \cite{boa} (1965) contained alternative proofs and further comments.
The authors were unaware of the earlier results.

Then came the work of Gromov \cite{cgt}, Anderson, Vamanamurthy, and Vuorinen
\cite{avv} \cite{avv2}, Pinelis \cite{pi1}, and many subsequent authors who have
done a tremendous amount of good work to enrich the subject area. 
The readers should have no difficulty finding them by searching for 
``monotone L'H\^opital rule'' on the internet, and by referring to the references
cited in the papers listed. 

Analogous to the classical L'H\^opital's rule for finding limits of indeterminant
forms, one can also consider the situation when $ f(0)=g(0)=\pm \infty  $. Results,
examples and counterexamples in this regard are presented
in Anderson, et al. \cite{avv2}. In this article,
we confine ourselves to the case when $ f(0) $ and $ g(0) $ are finite.

\section{Different formulations of the rules.}

If we let $ p(x)=f'(x) $ and $ q(x)=g'(x) $, then the rule stated in Section~1
becomes:

\par\vspace*{-3mm}\par
\begin{center}
\fbox{
\begin{minipage}[t]{5in}
\em 
\par\vspace*{0mm}\par
\CENTERLINE{If \,\, $\displaystyle \frac{p(x)}{q(x)} $ \,\, is monotone, so is \,\, $\displaystyle \frac{\vp{12}{6}\int_0^xp(t)\,dt}{\vp{12}{0}\int_0^xq(t)\,dt} $ ``in the same sense''.}
\par\vspace*{1mm}\par
\end{minipage}
}
\end{center}

This integral formulation has two advantages. 
One is that now we do not have to
separately assume that $ f(0)=g(0)=0 $. Second is that the rule, stated as
it is, can be applied to functions $ p $ and $ q $ with discontinuities.
For instance, $ p(x) $ and $ q(x) $ may be piecewise continuous, such as
step functions, as long as their quotient is still monotone.
In such cases, the differential form would have failed because the
requirement that $ f $ and $ g $ be differentiable at every point is not
satisfied. If one has proved the differential form of the rule previously based
on this assumption, such as invoking the generalized Cauchy mean value
theorem, one has to seek a different proof.

Alternatively, one can use an approximation technique, as suggested by 
Redheffer \cite{re}, to deduce the integral form from the differential form.
Choose a sequence of continuous functions $ q_n(x) $ that converge to $ q(x) $ in the uniform norm
in the finite interval $ [0,b\,] $. Next, choose another sequence of continuous functions $ h_n(x) $ that
converge to $ p(x)/q(x) $, this time requiring that each
$ h_n(x) $ preserves the same
monotone property of the latter. This can be done, for instance, by
using the well-known method of mollifiers of S. Sobolev and K.O. Friedrich
in the theory of partial differential equations. Since the
functions $ \vp{14}{8}\int_{0}^{x} h_n(t)q_n(t)\,dt $ and $ \int_{0}^{x} q_n(t)\,dt $ are now differentiable
everywhere, the differential form of the rule implies that
$ \vp{14}{8}\int_{0}^{x} h_n(t)q_n(t)\,dt /\int_{0}^{x} q_n(t)\,dt $ is monotone. 
Letting $ \vp{10}{0}n\rightarrow \infty  $ gives the desired conclusion.

There are other ways to relax the differentiability requirement
in the differential formulation. For example, one can use one-sided Dini
derivatives instead of regular derivatives. In the next section, we will
see how the rule can be formulated even without mentioning 
differentiability.
As we have remarked before, such extensions may be of theoretical value, but 
they are often not needed for practical applications.

Another direction of extension is to strengthen the conclusion of the
rule. Pinelis \cite{pi4} showed that if $ f'/g' $ is $ \ic $, then
$ (f/g)' $ is in fact strictly + (with strict
monotonicity following as a corollary).
One can compare this with the regular and strong forms of the maximum
principle in the theory of differential equations. The weaker form states that the global maximum
must be attained at the boundary of the region in consideration, while
the strong form maintains that at the boundary point where the global
maximum is attained, the directional derivative along the outward normal
must be strictly +.

One is also able to deduce strict monotonicity in a certain sense
from non-strict monotone assumptions. One such case will be discussed 
in the next section.

\section{Convex functions}

In the particular case when $ g(x)=x $, the increasing 
monotone rule reduces to:
\begin{quote}
\em If $ f(0)=0 $ and $ f'\nearrow $, so is $ f/x $.
\end{quote}
The slightly weaker (since $ f $ is required to be twice differentiable)
version, $ f''>0\,\,\Longrightarrow \,\,f(x)/x\ia $, is what
HLP referred to as Proposition 127. Both of these results are subsumed by
the following well-known property of convex functions.
\[
\fbox{
\begin{minipage}[t]{5in}
\par\vspace*{2mm}\par
\em If $ f $ is a continuous, strictly
convex function in $ [0,b) $, satisfying $ f(0)\leq 0 $,
then $ f/x\ia $ in $ [0,b) $.
\par\vspace*{1mm}\par
\end{minipage}
} 
\]

A function is defined to be convex  in $ [0,b\,] $ if 
for all $ 0\leq x_1<x_2\leq b $,
\begin{equation}  f\left( \frac{x_1+x_2}{2} \right)  \leq  \frac{f(x_1)+f(x_2)}{2} \,.  \end{equation}
It is said to be strictly convex if $ \leq  $
is replaced by $ < $. Concavity is defined by reversing the inequality signs.
No differentiability requirement is assumed. 

It is well-known that
for a continuous convex function, the following inequality holds:
\begin{equation}  f(\lambda x_1+(1-\lambda )x_2) \leq  \lambda  f(x_1) + (1-\lambda ) f(x_2), \qquad  0<\lambda <1 .  \Label{con}  \end{equation}
If the convexity is strict, replace $ \leq  $ by $ < $.
Its geometric interpretation is that the arc of the
graph of $ f(x) $, $ x\in[x_1,x_2] $, lies below the the chord joining the
two points $ A_1=(x_1,f(x_1)) $ and $ A_2=(x_2,f(x_2)) $.

\begin{center}
\FG{50}{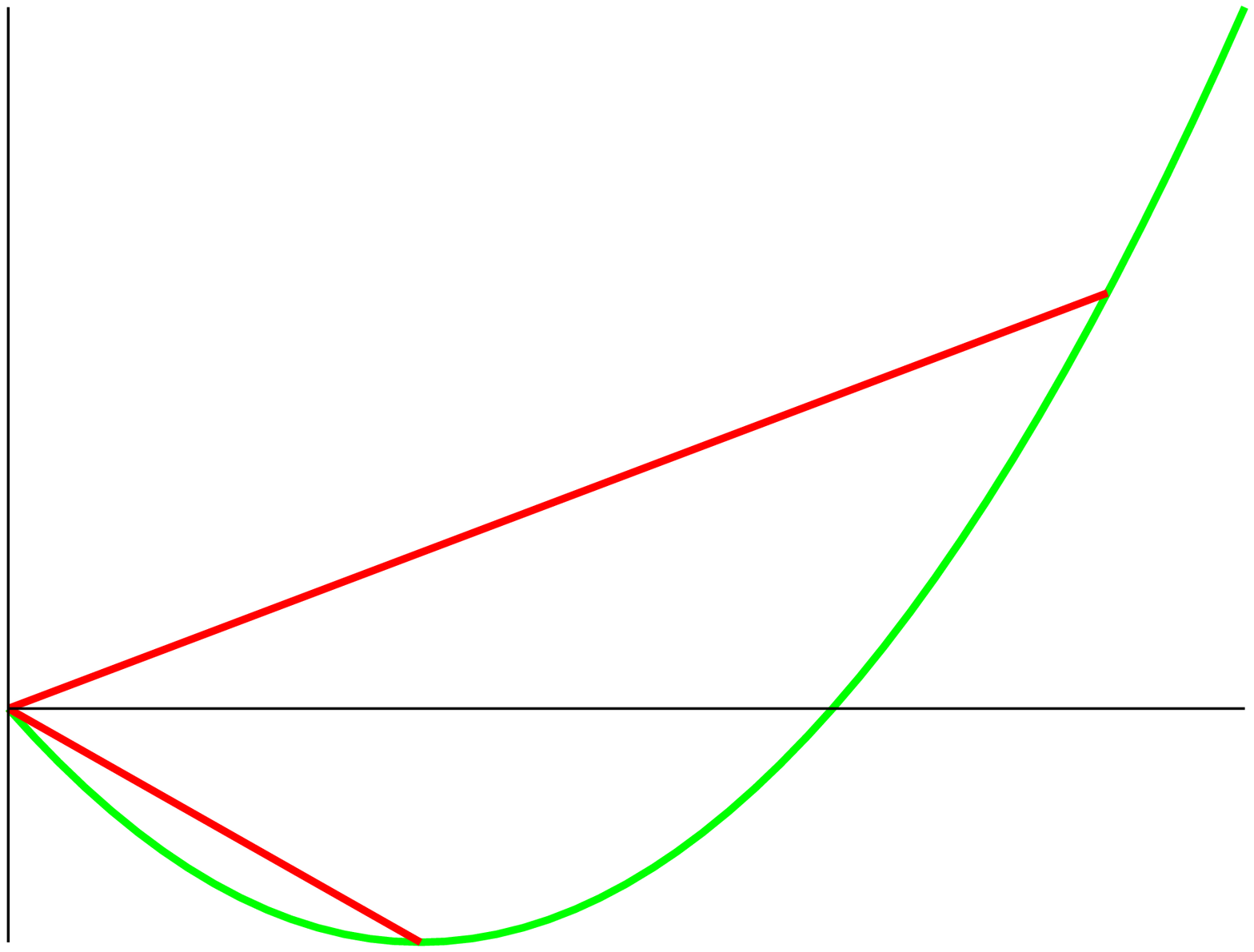} \qquad  \qquad  \FG{50}{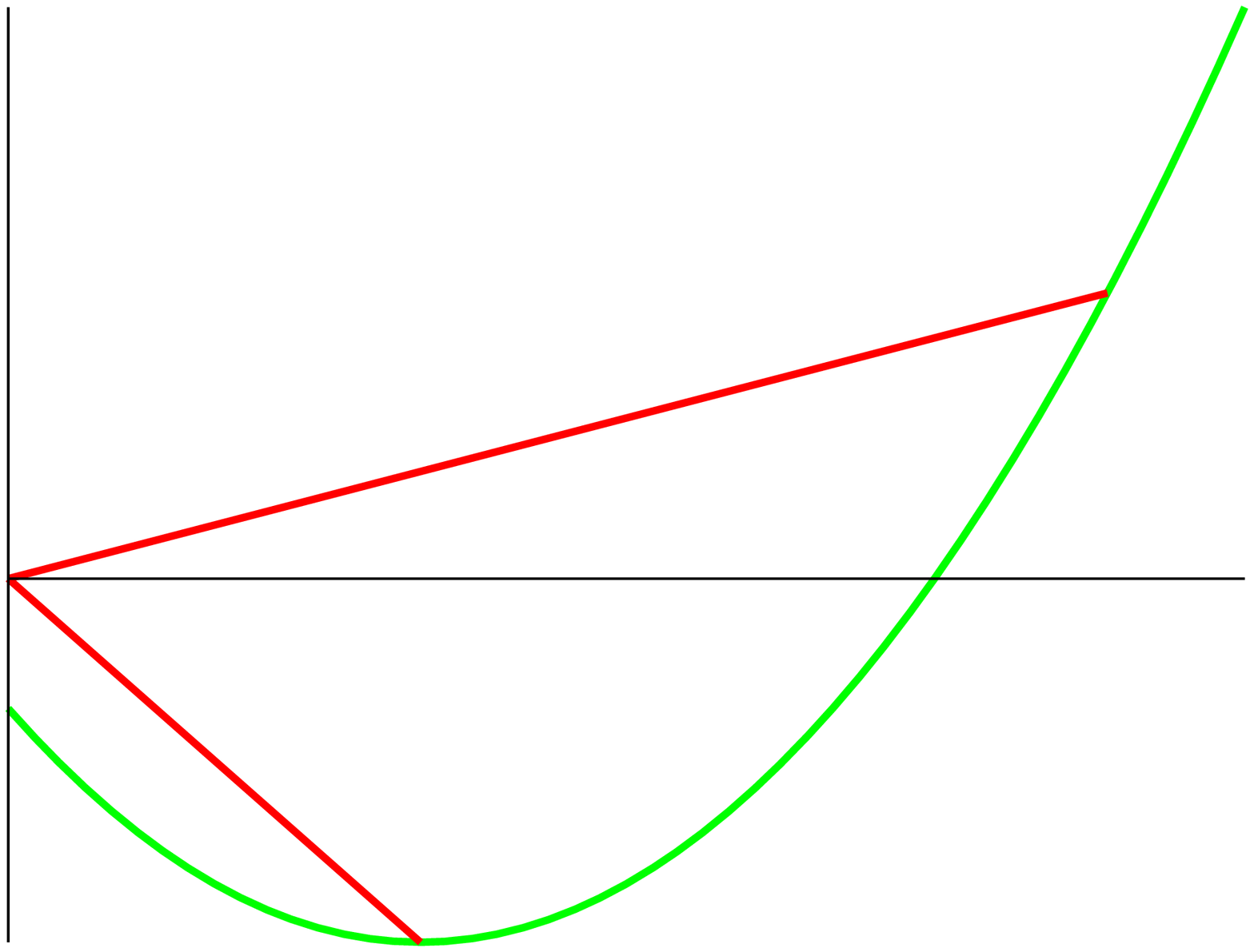}

\par\vspace*{6mm}\par
Figure 1. Convex functions with $ f(0)=0 $ and $ f(0)<0 $.
\end{center}

\par\vspace*{-31mm}\par
\hspace*{7mm} $ O $
\par\vspace*{-14mm}\par
\hspace*{75mm} $ O $
\par\vspace*{9mm}\par
\hspace*{32mm} $ A_1 $
\par\vspace*{-6mm}\par
\hspace*{98mm} $ A_1 $
\par\vspace*{-44mm}\par
\hspace*{57mm} $ A_2 $
\par\vspace*{-7mm}\par
\hspace*{124mm} $ A_2 $
\par\vspace*{50mm}\par

Figure 1 depicts two such functions (green curves); 
the first has $ f(0)=0 $ and the second $ f(0)<0 $.
The quantity $ f(x)/x $ represents the slope of the straight line $ OA $
joining the origin $ O $ and the point $ A=(x,f(x)) $ on the curve.
As $ x $ increases, the point $ A $ slides along the curve towards the
right and it is intuitively clear that
the slope of $ OA $ increases. A vigorous proof can be
given using (\ref{con}). The proof is not new, but we include it here for
easy reference.

\begin{PROOF}
In the case $ f(0)=0 $, by convexity, 
the arc of the curve between $ O $ and $ A_2 $
lies below the straight line $ OA_2 $. In particular, the point $ A_1 $ lies
below $ OA_2 $ and the desired conclusion follows.

For the case $ f(0)<0 $, we 
modify the function $ f $ in $ [0,x_1] $ by replacing the arc 
over $ [0,x_1] $ with
the chord $ 0A_1 $. The resulting new curve represents the function
$ \max(f(x),f(x_1)x/x_1) $ which is again convex.
Then we are back to the first case $ f(0)=0 $.
\end{PROOF}

{\bf Alternative Proof}.
Another simple proof uses the fact that a straight line
cannot intersect a strict convex{/}concave curve at more than two points.
As above, we only have to consider the case $ f(0)=0 $.
Suppose $ f/x $ is not monotone; then there are $ x_1\neq x_2 $ such that
$ f(x_1)/x_1=f(x_2)/x_2=\lambda  $. Then the straight line $ y=\lambda x $ intersects
the graph of $ f $ at the two points $ (x_1,f(x_1)) $ and $ (x_2,f(x_2)) $.
A third intersection point, however, is $ (0,0) $, giving a contradiction.
{{\ \vrule height7pt width4pt depth1pt} \par \vspace{2ex} }

But how do we make the quantum jump from the special case $ g(x)=x $
to the general case?
The trick of change of variable, just as HLP pointed out in their proof, 
is the key. By the way, this trick has been known to work for the
classical L'H\^opital rule for indeterminant limits as well. 
See, for example, Taylor \cite{ta}. Let us rephrase
HLP's argument to make it more transparent. 

Suppose that $ g(0)=0 $ and $ g(x) $ is continuous and $ \ic $ in $ [0,b\,] $.
The inverse function $ x=g^{-1}(u) $ is well-defined.
Substitute this into the definition of $ f(x) $ to get the composite function
$ \phi (u)=f(g^{-1}(u)) $. We have then the following generalized monotone rule.

\par\vspace*{-3mm}\par
\begin{center}
\fbox{
\begin{minipage}[t]{5in}
\par\vspace*{2mm}\par
\em If $ \phi (u) $ is a continuous, strictly convex {\rm(}concave{\rm)} function of $ u $ in
$ [0,g(b)] $ and $ f(0)\leq 0 $ $ (\geq 0) $,
then $ \displaystyle \vp{20}{0}\frac{f(x)}{g(x)} =\frac{\phi (u)}{u} $ is $ \ic $ {\rm(}$\dc${\rm)}.
\par\vspace*{2mm}\par
\end{minipage}
}
\end{center}

In practical
applications, to check the convexity or concavity of $ \phi (u) $, we 
often resort to showing that 
$ \displaystyle \phi '(u)  $
is monotone. By the chain rule, 
$ \displaystyle \phi '(u)  $
is nothing but $ \displaystyle \frac{f'(x)}{g'(x)} $.

This new rule is more general than the popular
one described on p.\ \!\!1,
because it does not impose any differentiability on $ f $ and $ g $, and
it covers situations when $ f(0)\neq 0 $.

If we omit ``strictly'' in the hypotheses, we cannot guarantee that $ f/g $ is
strictly monotone, as the degenerate example $ f(x)=g(x)=x $ shows.
However, this is pretty much the only exceptional situation, in the following
sense. If there exist two points $ x_1<x_2 $ in $ [0,b) $ such that 
$ f(x_1)/g(x_1)=f(x_2)/g(x_2) $, then $ f(x)=\lambda g(x) $ in $ [0,x_2] $
for some constant $ \lambda  $. To prove this, one only has to study the special case
when $ g(x)=x $. Then the assertion becomes geometrically obvious. We omit 
the details.

\begin{center}
\FG{50}{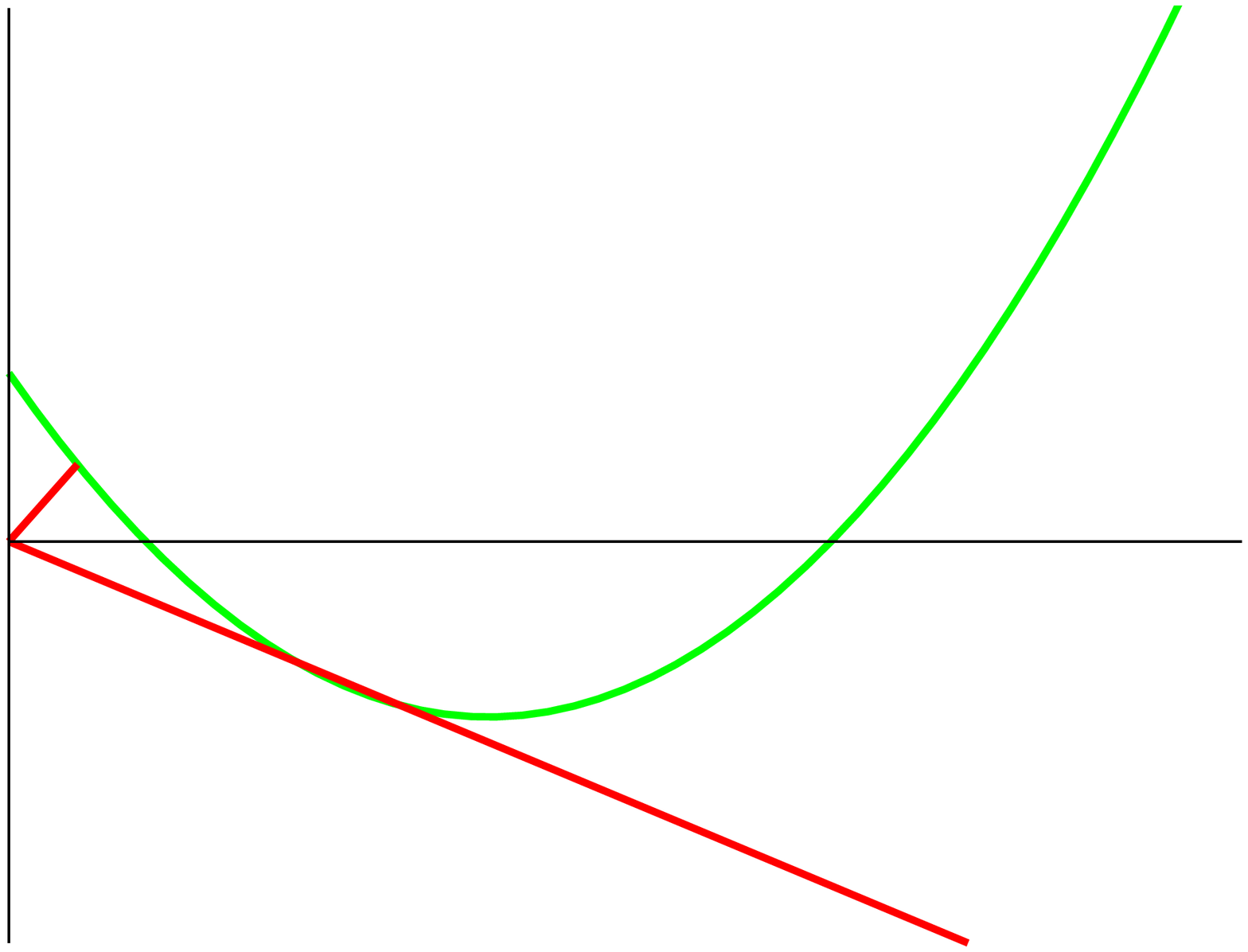}

\par\vspace*{2mm}\par
Figure 2. Convex function with $ f(0)>0 $.
\end{center}

\par\vspace*{-40mm}\par
\hspace*{32mm} $ O $
\par\vspace*{-14mm}\par
\hspace*{42mm} $ A $
\par\vspace*{12mm}\par
\hspace*{52mm} $ A_c $

\newpage
Next, let us look at the case  $ f(0)>0 $ when $ g(x)=x $, and $ f(x) $ is strictly
convex. Figure 2 depicts such a curve.

Starting at $ x=0 $, the point $ A $ is where the curve intersects the vertical axis
and the slope of $ OA $ is $ \infty  $. As $ x $ increases, the slope of $ OA $ decreases 
until we reach the point $ A_c $ such that $ OA_c $ is tangent to the curve.
Beyond that the slope of $ OA $ increases again. This is an example of 
the general situation described in the conclusion of HLP's Proposition.
The quintessential conclusion is that there exists at most one point $ c\in(0,b) $
such that $ f/x $ is $ \dc $ in $ (0,c) $ and $ \ic $ in $ (c,b) $.
As the example $ f(x)=1/(1+x) $ shows, the point $ c $ may not exist at all. 
After uplifting to the case of general $ g(x) $, the pertinent part of
HLP's Proposition can be generalized as follows. The symbol $ \exists\,!  $
is a shorthand for  ``there exists a unique''.

\par\vspace*{-3mm}\par
\begin{center}
\fbox{
\begin{minipage}[t]{5in}
\par\vspace*{2mm}\par
\em $ \phi (u) $ is continuous, strictly convex {\rm(}concave{\rm)} in $ u $ and $ f(0)>0 $ $ (<0) $.

Then either $ \displaystyle \vp{20}{0}\frac{f(x)}{g(x)} $ is $ \dc $ {\rm(}$\ic${\rm)} in $ (0,b) $, or 

$ \exists\,!  $
$ c\in(0,b) $ such that $ \displaystyle \vp{20}{0}\frac{f(x)}{g(x)} $ is $ \dc $ {\rm(}$\ic${\rm)} in $ (0,c) $ and $ \ic $ {\rm(}$\dc${\rm)} in $ (c,b) $.
\par\vspace*{2mm}\par
\end{minipage}
}
\end{center}

\begin{PROOF}
Of course, HLP's proof no
longer works under the minimal assumption. A rigorous proof
can be given using only properties of convex functions.
The quotient $ f/x $ is continuous in $ (0,b\,] $. Although
it blows up at $ x=0 $, it is easy to see that
it attains a global minimum at some point $ c\in(0,b\,] $.
By its very construction,
the line $ OA_c $ that joins the origin $ O $ and the point $ A_c=(c,f(c)) $
lies below the curve of $ f $. We can show that $ f/x $ is decreasing in $ [0,c] $
as follows. Let $ 0<x_1<x_2<c $.
Since the line joining $ O $ and the point $ A_1=(x_1,f(x_1)) $
is above the line $ OA_c $, the former must intersect the vertical line $ x=c $
above $ A_c $, say at a point $ D $.
The arc of $ f(x) $ from $ A_1 $ to $ A_c $
must be below the line $ A_1A_c $, which lies below $ A_1D $. In particular, the 
point $ A_2=(x_2,f(x_2)) $ is below the
line $ A_1D $ and hence $ f(x_2)/x_2<f(x_1)/x_1 $.

If it happens that $ c=b $,
$ f/x $ has no chance to bounce back from
$ \dc $ to $ \ic $. If $ c<b $, then we can argue in a similar way as above that
$ f/x $ is $ \ic $ in $ (c,b) $. This completes the proof.
\end{PROOF}

Is there a practical way to delineate the two situations in the conclusion?
In the special case $ g(x)=x $ and $ f(x) $ is convex, it is easy to see that
a necessary and sufficient condition for the existence of an interior $ c\in(0,b) $ is that
in a left neighborhood of $ b $, the curve of $ f $ lies below the chord joining 
$ O $ and $ B=(b,f(b)) $. In the general setting, this criteria can be expressed as
\begin{equation}  \liminf _{x\rightarrow b^-} \, \frac{f(b)-f(x)}{g(b)-g(x)} > \frac{f(b)}{g(b)} \,.  \end{equation}
One can also use limsup instead of liminf.
For differentiable $ f(x) $ and $ g(x) $, it simplifies to
$$  f(b) g'(b) - f'(b) g (b) < 0 \, .  $$
For concave $ \phi (u) $, the inequality sign is reversed. If $ c<b $ exists, it is 
determined by 
solving the equation $ f'(x)g(x)=f(x)g'(x) $.

If we allow the function $ \phi (u) $ to be non-strict convex in the hypotheses, then 
the curve of $ \phi (u) $ (versus $u$) may contain flat portions 
(line segments). If it
happens that the tangent line from the origin touches the curve
and contains one of these flat portions, 
then the unique turning point $ c $ in the above
rule becomes an entire interval of turning points. We have to modify the rule to 
say that now there exists a subinterval $ [c_1,c_2]\subset(a,b) $ such that
$ f/g $ is $ \dc $ in $ (0,c_1) $, constant in $ [c_1,c_2] $ and $ \ic $ in $ (c_2,b) $.

Finally, let us consider the case when $ g(0)>0 $, which is also covered
by HLP's Proposition. In the special case when
$ g(x)=x+\gamma  $, $ \gamma >0 $, and $ f $ is convex,
the corresponding problem is to investigate the $ \ic $ and $ \dc $
properties of $ f(x)/(x+\gamma ) $, which is the slope of the line joining the
point $ (-\gamma ,0) $ on the $ x $-axis and the point $ (x,f(x)) $ on the graph of $ f $.
Equivalently, we may apply a translation to shift the point $ (-\gamma ,0) $
to be the new origin. In this perspective, we can exploit the same figures
earlier in this section, only that now the curve of $ f $ starts from $ x=\gamma  $
instead of $ x=0 $.

It is easy to check with simple examples that both
possibilities discussed in the case $ g(0)=0 $, $ f(0)>0 $ can occur. 
Besides those, an additional
possible third situation is that $ f/g $
is $ \ic $ in $ (0,b) $. In the special case $ g(x)=x+\gamma  $, this happens when
the line joining the points $ (-\gamma ,0) $ and $ (0,f(0)) $ lies below the
curve of $ f $ in a right neighborhood of $ x=0 $. In the general setting, if $ f $ and
$ g $ are differentiable, this happens when
$ f(0)/g(0) - f'(0)/g'(0) \leq  0. $
However, the possibility of $ f/g $ having the shape $ \ic\dc $ is not allowed.

We summarize all the findings:

\begin{THEOREM}
Suppose that $ g(x) $ is a continuous, positive $ \ic $ function in $ [0,b\,] $ and $ \phi (u)=f(g^{-1}(u)) $
is a continuous strictly convex {\rm(}concave{\rm)} function of $ u $ in $ [g(0),g(b)] $.

\par\vspace*{-2mm}\par
\begin{itemize}
\item[\rm(1)\,\,] $ g(0)=0 $ and $ f(0)\leq (\geq )\,\,0 $,
   then $ f/g $ is $ \ic $ $(\dc)$.
\item[\rm(2)\,\,] $ g(0)=0 $ and $ f(0)>(<)\,\,0 $\quad  or \quad  $ g(0)>0 $.
   \begin{itemize}
   \item[\rm i)\,\,] If \,\, $ g(0)\neq 0 $ and 
$ \displaystyle \left(  \frac{f(0)}{g(0)} - \liminf_{x\rightarrow 0+} \, \frac{f(x)-f(0)}{g(x)-g(0)} \right)  \leq  (\geq 0)\,\, 0 , $
       then $ f/g $ is $ \ic $ $ (\dc) $.
   \item[\rm ii)\,\,] If \,\,
$ \displaystyle \left(  \frac{f(b)}{g(b)} - \liminf_{x\rightarrow b-} \, \frac{f(b)-f(x)}{g(b)-g(x)} \right)  \geq  (\leq 0)\,\, 0 , $
       then $ f/g $ is $ \dc $ $ (\ic) $.
   \item[\rm iii)\,\,] Otherwise, 
       $ f/g $ has the shape $ \dc\ic $ $ (\ic\dc) $ with a unique turning point in $ (0,b) $.
   \end{itemize}
\end{itemize}

Suppose that the convexity{/}concavity of $ \phi (u) $ is not assumed to be 
strict. Let $ [0,\alpha )\subset[0,b\,] $ be a maximal 
subinterval $($which is possibly void$)$
in which $ f(x)=\lambda g(x) $ for a constant $ \lambda  $. Then in $ [\alpha ,b\,] $, the same conclusions
as above $($with all monotonicity being strict$)$ hold.
\end{THEOREM}

\begin{REM} \rm
Case (1) can actually be combined with (2) i). We prefer to separate it
out since it is historically as well as application-wise 
the most prominent case.
\end{REM}

\begin{REM} \rm
Two points in the result are worth noting. The function $ f/g $ cannot change monotonicity
more than once, and in the convex case, the shape $ \ic\dc $ is ruled out.
\end{REM}

\begin{REM} \rm
The usefulness of the criteria given in the cases (2) i) and ii) 
(in addition to the $ f'/g'\ic $ condition) lies 
in the fact that there is only one boundary condition at
one of the endpoints to verify in order to 
deduce monotonicity of $ f/g $ over the entire interval.
\end{REM}

\begin{REM} \rm
The last part of the Theorem is a strong form of the rule. Even if strict
convexity{/}concavity is not assumed, we can still obtain
strict monotonicity, except possibly in an initial subinterval in a very
special situation.
\end{REM}

If $ g(x) $ is a positive $ \dc $ function in $ (0,b) $, an analogous result holds. We only have
to do a reflection $ x\mapsto (b-x) $ to reduce it to the $ \ic $ case. The
role of $ 0 $ and $ b $ are now exchanged. One has to be careful in chasing
the signs and inequalities in the conditions. We state the result for ease of
reference. It is useful in handling functions such as $ f(x)/(b-x) $.

\begin{THEOREM}
Suppose that $ g(x) $ is a continuous, positive $ \dc $  function in $ (0,b) $ and $ \phi (u)=f(g^{-1}(u)) $
is a continuous strictly convex $($concave$)$ function of $ u $ in $ [g(b),g(0)] $.

\par\vspace*{-2mm}\par
\begin{itemize}
\item[\rm(1)\,\,] $ g(b)=0 $ and $ f(b)\leq (\geq )\,\,0 $,
   then $ f/g $ is $ \dc $ $(\ic)$.
\item[\rm(2)\,\,] $ g(b)=0 $ and $ f(b)>(<)\,\,0 $\quad  or \quad  $ g(b)>0 $.
   \begin{itemize}
   \item[\rm i)\,\,] If \,\,
$ \displaystyle \left(  \frac{f(0)}{g(0)} - \liminf_{x\rightarrow 0+} \, \frac{f(x)-f(0)}{g(x)-g(0)} \right)  \geq  (\leq 0)\,\, 0 , $
       then $ f/g $ is $ \ic $ $ (\dc) $.
   \item[\rm ii)\,\,] If \,\, $ g(b)\neq 0 $\,\, and
$ \displaystyle \left(  \frac{f(b)}{g(b)} - \liminf_{x\rightarrow b-} \, \frac{f(b)-f(x)}{g(b)-g(x)} \right)  \leq  (\geq 0)\,\, 0 , $
       then $ f/g $ is $ \dc $ $ (\ic) $.
   \item[\rm iii)\,\,] Otherwise, 
       $ f/g $ has the shape $ \dc\ic $ $ (\ic\dc) $ with a unique turning point in $ (0,b) $.
   \end{itemize}
\end{itemize}
\end{THEOREM}

In practical applications when $ f $ and $ g $ are differentiable, the
following simplified rule is easier to use. Corollary 1 deals with increasing
$ g $ and Corollary 2 deals with decreasing $ g $.

\begin{COR}
Suppose $ f $ and $ g $ are differentiable,
$ g,\,g'>0 $, and $ f'/g'\ic(\dc) $ in $ (0,b) $.

\par\vspace*{-2mm}\par
\begin{itemize}
\item[\rm(1)\,\,] If $ g(0)=0 $ and $ f(0)\leq (\geq )\,\,0 $, then $ f/g $ is $ \ic $ $(\dc)$.
\item[\rm(2)\,\,] $ g(0)=0 $ and $ f(0)>(<)\,\,0 $\quad  or \quad  $ g(0)>0 $.
\par\vspace*{-2mm}\par
   \begin{itemize}
   \item[\rm i)\,\,] If \,\, $ g(0)\neq 0 $\,\, and
$ \displaystyle (f/g)(0) \leq  (\,\geq \,)\,  (f'/g')(0), $
       then $ f/g $ is $ \ic $ $ (\dc) $.
   \item[\rm ii)\,\,] If \,\,
$ \displaystyle (f/g)(b) \geq  (\,\leq \,)\, \hspace*{0.6mm} (f'/g')(b), $
       then $ f/g $ is $ \dc $ $ (\ic) $.
   \item[\rm iii)\,\,] Otherwise, 
       $ f/g $ has the shape $ \dc\ic $ $ (\ic\dc) $ with a unique turning point in $ (0,b) $.
   \end{itemize}
\end{itemize}
\end{COR}

\newpage
\begin{REM} \rm
It is instructive to visualize the three possibilities in the rule.
Figure 3 shows the plots of $ f/g $ (red curves) and $ f'/g' $
(green dashed curves)
in three typical examples. All other examples exhibit the same features. 
By hypotheses, the dashed curves are $ \ic $.
\end{REM}

\begin{center}
\FG{42}{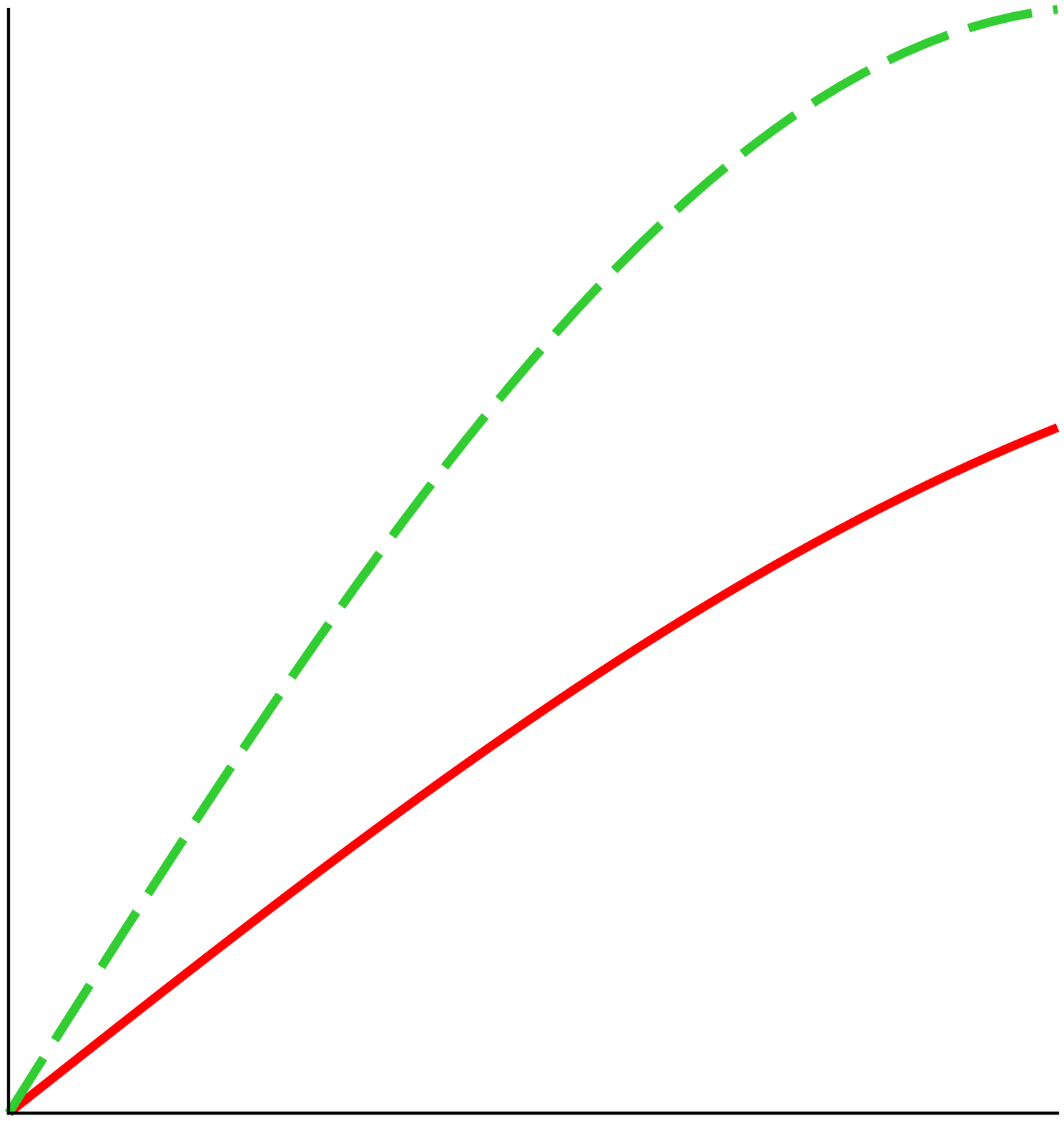} \qquad  \FG{42}{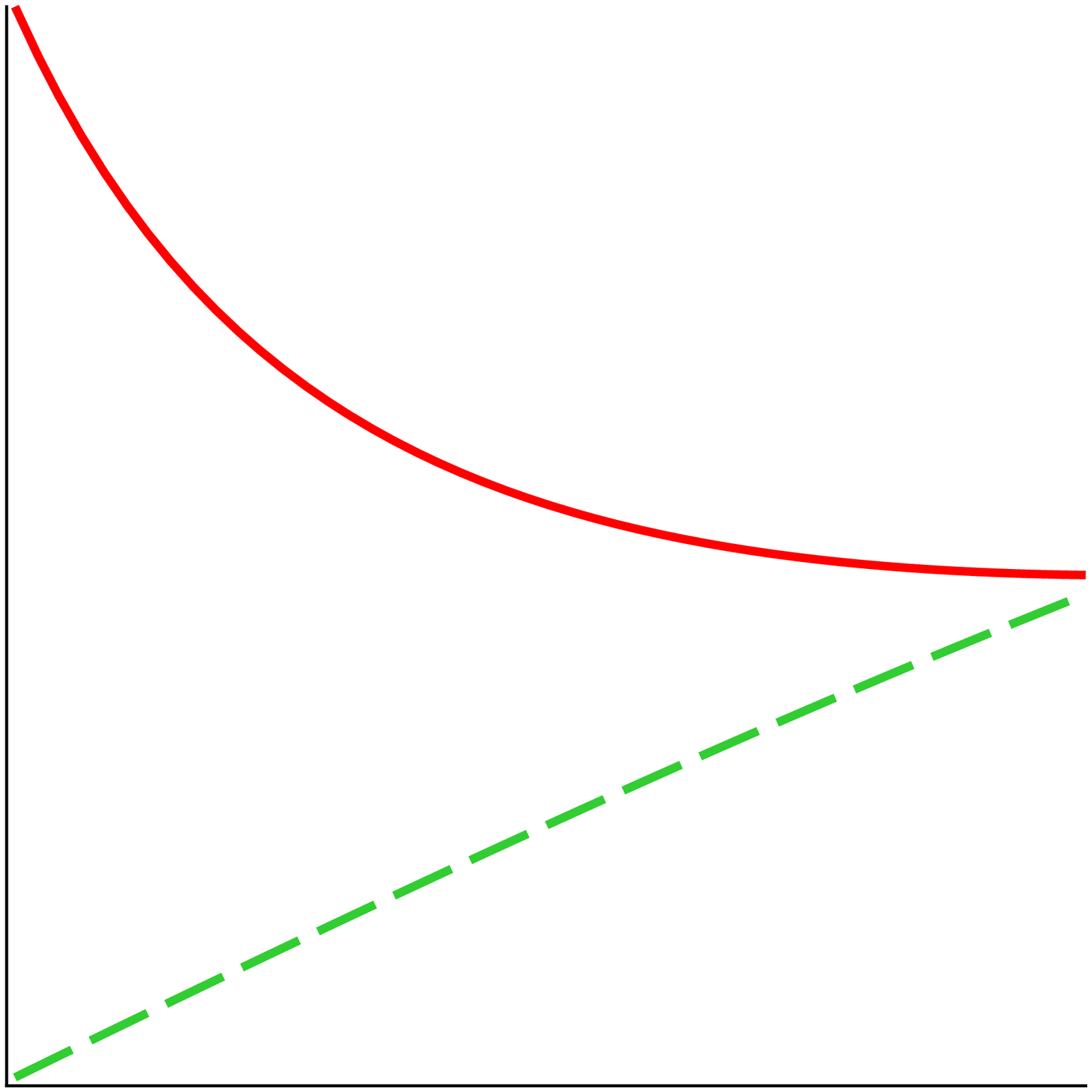} \qquad  \FG{42}{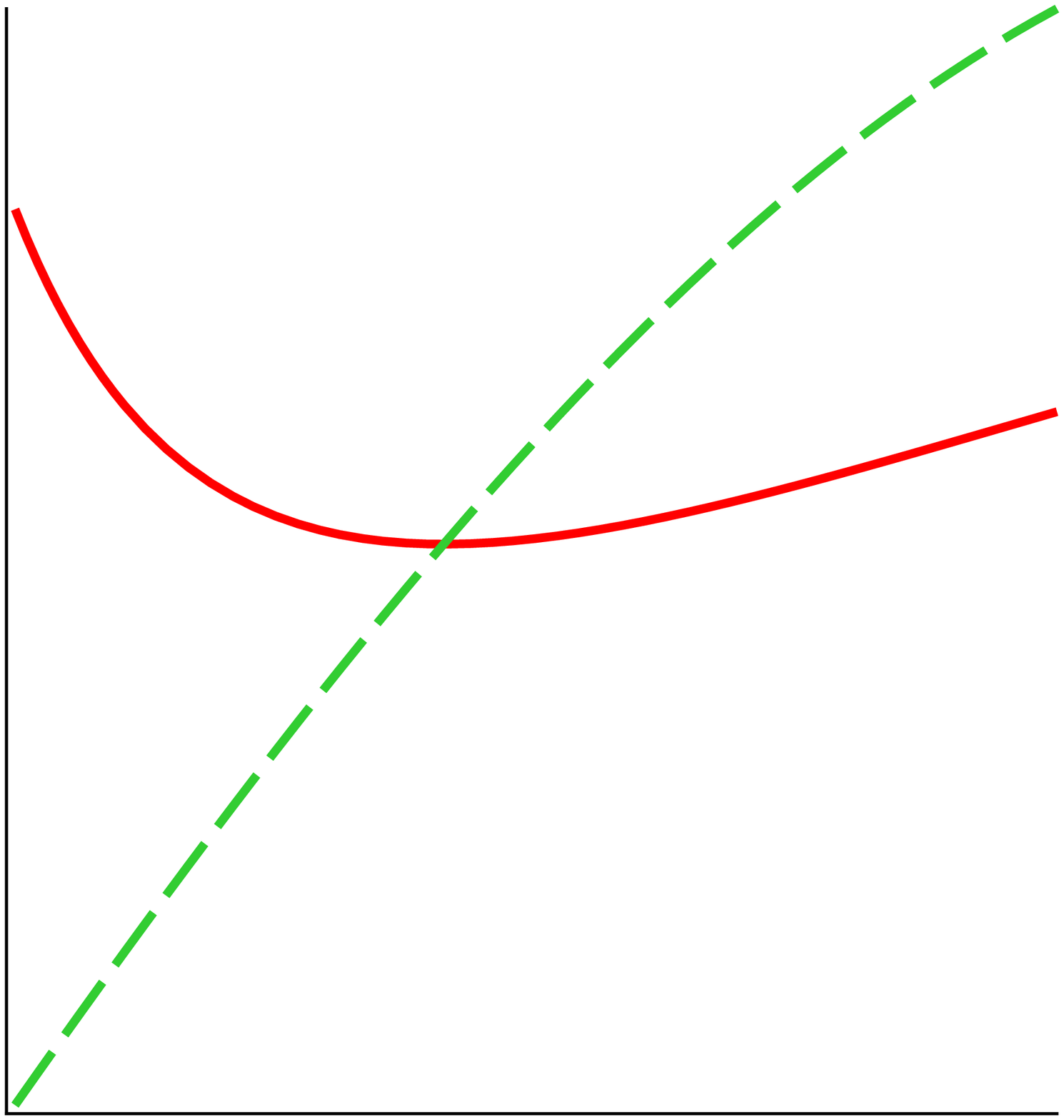}

 Cases (1) and (2) i) \hspace*{26mm} (2) ii) \hspace*{40mm} (2) iii)  \hspace*{5mm}

\par\vspace*{6mm}\par
Figure 3. Three possibilities of $ f/g $ (red curve) when $ f'/g'\ia $.
\end{center}

\par\vspace*{1mm}\par
The first plot shows that if $ f/g $ lies below $ f'/g' $, then $ f/g $ is $ \ic $.
As the rule asserts, to guarantee this situation, you only have to check
whether $ f(0)=g(0)=0 $ (case (1)) or in case (2) i), to check whether
$ f(0)/g(0)\leq f'(0)/g'(0) $. In other words, only the behaviors at the left
endpoint matters.

The second plot shows that if $ f/g $ lies above $ f'/g' $, then $ f/g $ is $ \dc $.
Again, to guarantee this situation, you only need to 
check whether $ f(b)/g(b)\geq f'(b)/g'(b) $ at the right endpoint.

The third plot shows the $ \dc\ic $ possibility, which happens when the dashed
curve intersects the red curve at some point $ c\in(0,b) $. This situation can
be considered as a hybrid case: before the intersection, the plot looks
like the second one, and after that it looks like the first one.
Right at the intersection, the red curve has a horizontal tangent.
HLP have recorded the same observation in their proof of Proposition 148.
\par\vspace*{3mm}\par
\begin{COR}
Suppose $ f $ and $ g $ are differentiable, 
$ g>0 $, $ g'<0 $, and $ f'/g'\dc(\ic) $ in $ (0,b) $.

\par\vspace*{-2mm}\par
\begin{itemize}
\item[\rm(1)\,\,] If $ g(b)=0 $ and $ f(b)\leq (\geq )\,\,0 $, then $ f/g $ is $ \dc $ $(\ic)$.
\item[\rm(2)\,\,] $ g(b)=0 $ and $ f(b)>(<)\,\,0 $\quad  or \quad  $ g(b)>0 $.
\par\vspace*{-2mm}\par
   \begin{itemize}
   \item[\rm i)\,\,] If \,\,
$ \displaystyle (f/g)(0) \geq  (\,\leq \,)\,  (f'/g')(0), $
       then $ f/g $ is $ \ic $ $ (\dc) $.
   \item[\rm ii)\,\,] If \,\, $ g(b)\neq 0 $\,\, and
$ \displaystyle (f/g)(b) \leq  (\,\geq \,)\, \hspace*{0.6mm} (f'/g')(b), $
       then $ f/g $ is $ \dc $ $ (\ic) $.
   \item[\rm iii)\,\,] Otherwise, 
       $ f/g $ has the shape $ \dc\ic $ $ (\ic\dc) $ with a unique turning point in $ (0,b) $.
   \end{itemize}
\end{itemize}
\end{COR}

\newpage
\begin{REM} \rm
The analogous plots of $ f/g $ and $ f'/g' $ are shown in Figure 4.
\end{REM}

\begin{center}
\FG{42}{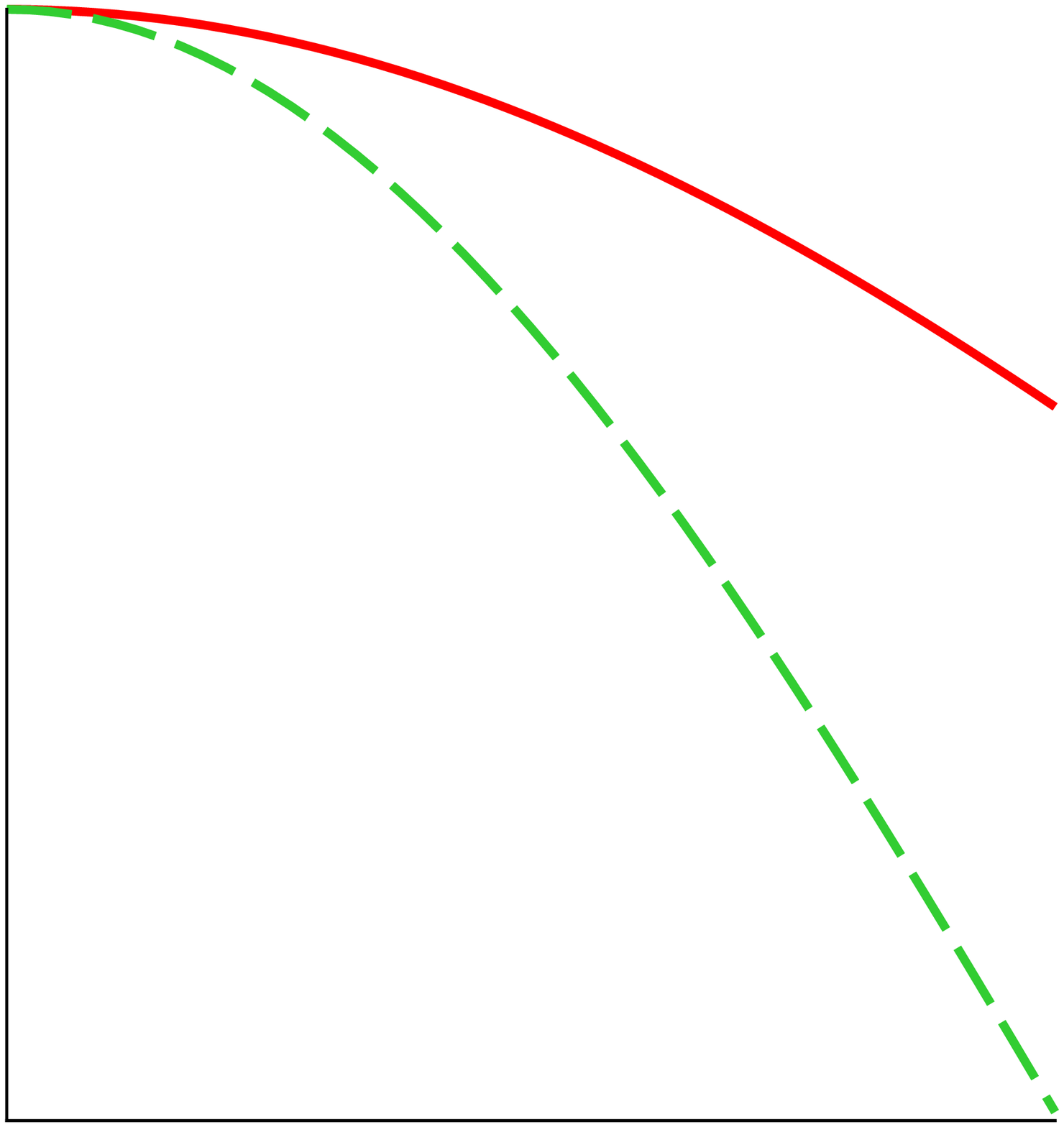} \qquad  \FG{42}{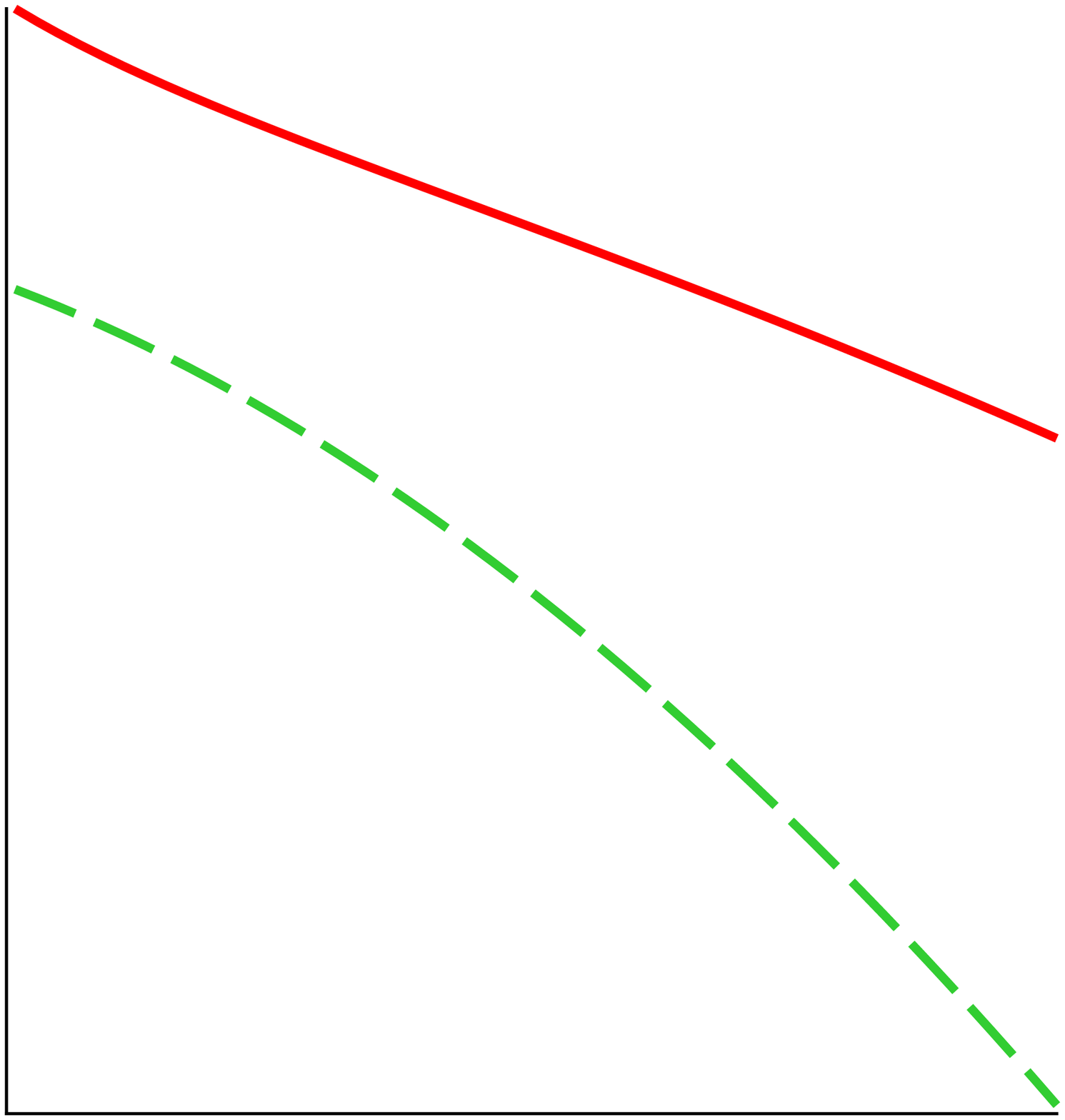} \qquad  \FG{42}{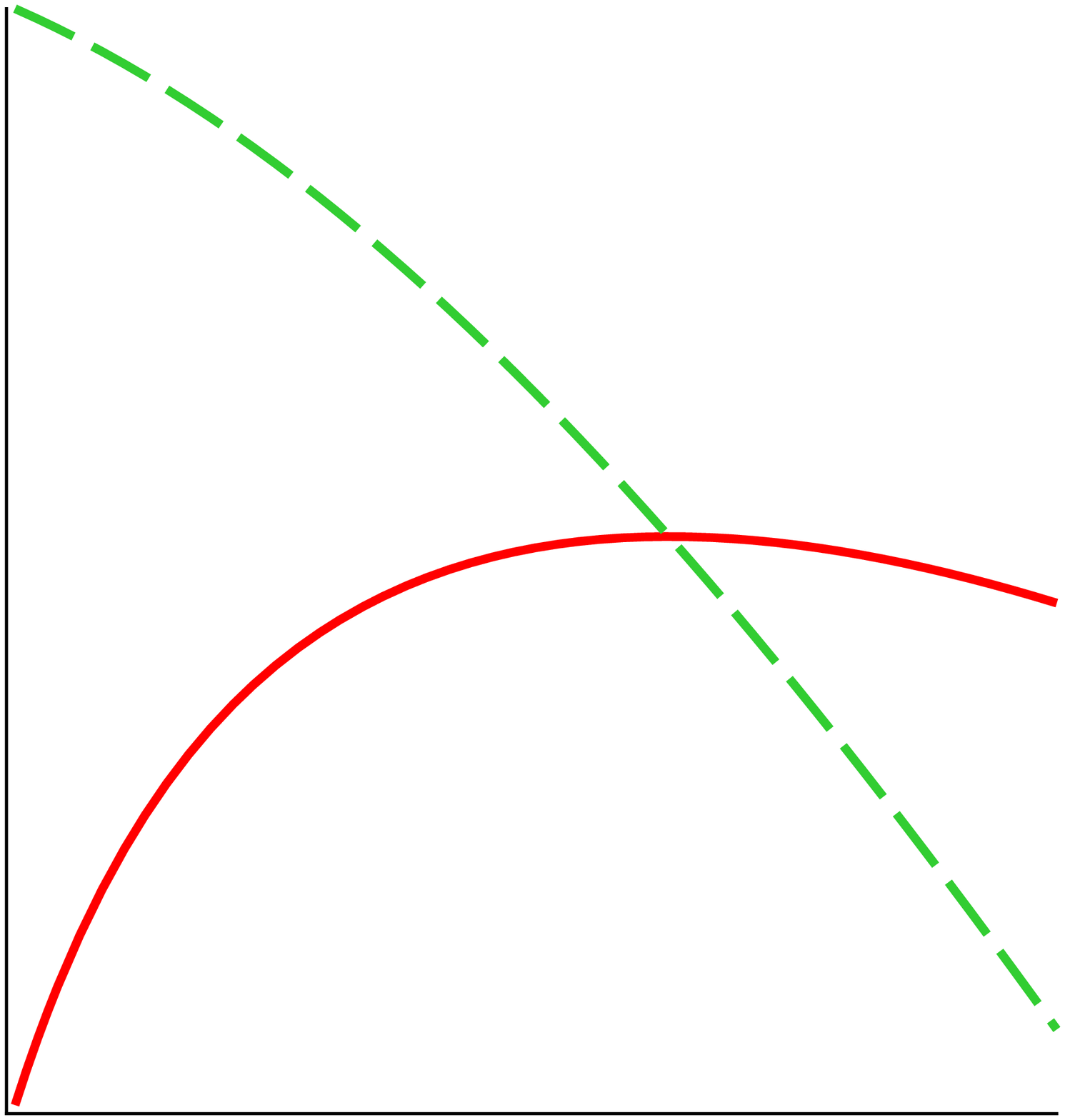}

 Cases (1) and (2) i) \hspace*{26mm} (2) ii) \hspace*{40mm} (2) iii)  \hspace*{5mm}

\par\vspace*{6mm}\par
Figure 4. Three possibilities of $ f/g $ (red curve) when $ f'/g'\da $.
\end{center}

\section{Rules for oscillation}

Pinelis \cite{pi3} \cite{pi4} discovered an interesting extension of the 
monotone rule to the situation when $ f'/g' $ is no longer
monotone. In this section, we look at this extension from the
perspective of convex{/}concave functions.

\begin{center}
\FG{50}{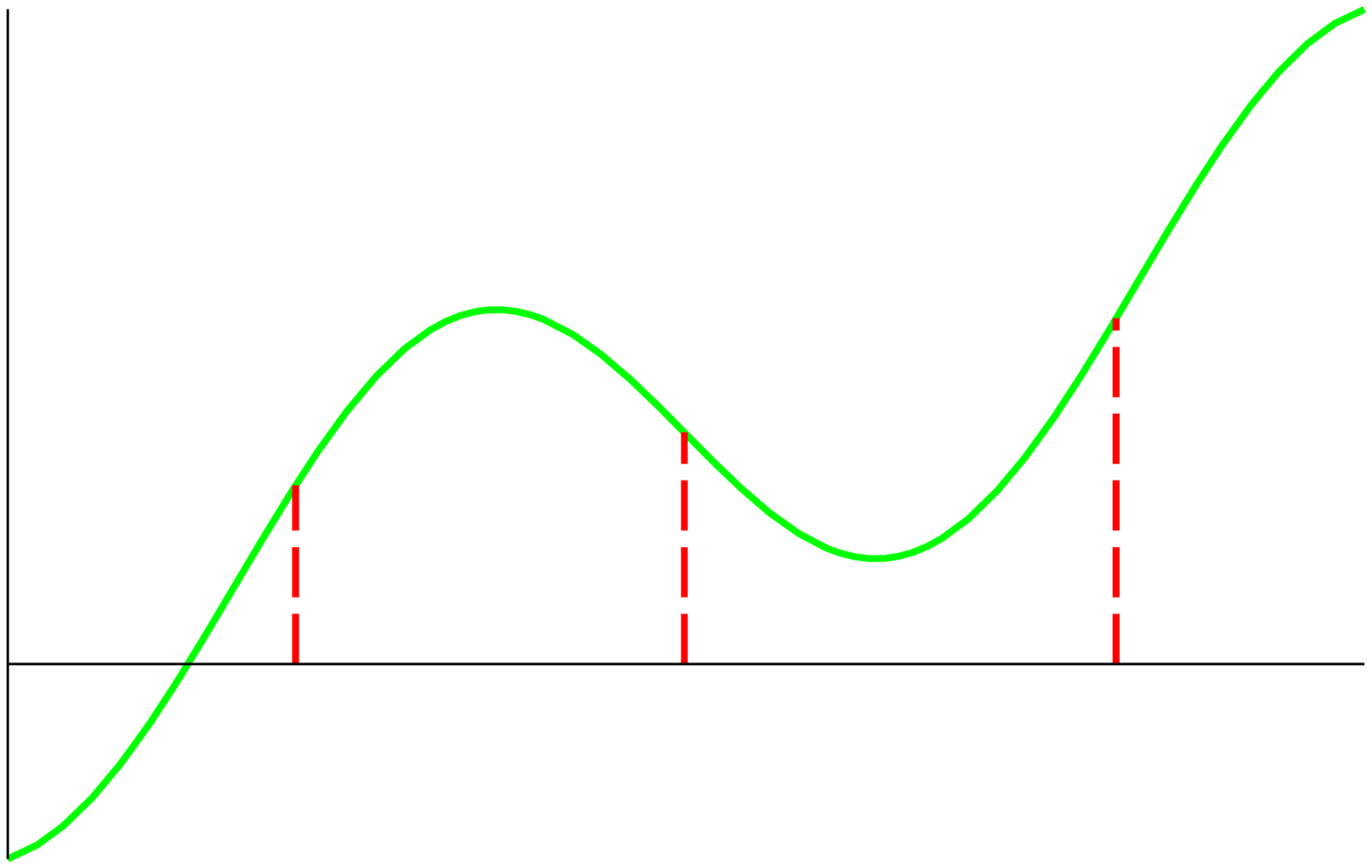}

Figure 5. A function that changes convexity three times.
\end{center}

\par\vspace*{-29mm}\par
\hspace*{24mm} $ O $
\par\vspace*{-3mm}\par
\hspace*{47mm} $ b_1 $\hspace*{18mm} $ b_2 $\hspace*{19mm} $ b_3 $
\par\vspace*{20mm}\par

Suppose that $ (0,b) $ can be divided into $ n+1 $ subintervals with points
$$  (0=b_0)\,<\,b_1\,<\,b_2\,<\,\cdots\,<\,b_n\,<\,(b_{n+1}=b),   $$
and $ f'/g' $ is assumed to be $ \ic $
in the odd subintervals $ (0,b_1) $, $ (b_2,b_3) $, $ \cdots $ and $ \dc $ in
the even ones. The function changes its monotonicity $ n $ times.
Figure 5 depicts one such functions with $ f(0)\leq 0 $. It is convex in 
$ [0,b_1] $, concave in $ [b_1,b_2] $, etc. The points $ b_i $ are points of reflection
of $ f $.

For convenience, we can also say that $ f'/g' $\, ``oscillates''
$ n $ times. Do not confuse this use of the term with the more
conventional meaning, as in the ``oscillation'' of the pendulum. When
a pendulum oscillates 2 times, it has ``oscillated in our sense''
(changed directions) 3 times.
The case when $ f'/g' $ is $ \dc $ in the odd subintervals and $ \ic $
in the even ones can be studied in a similar way.

\underline{$\vphantom{y}$Question}: {\em What can we say about the oscillatory property of
$ f/g $? }

Again, we first appeal to the special case when $ g(x)=x $,
and $ f $ is assumed to be alternatively convex and concave in the subintervals.
By Theorem~1 (1), $ f/x $ is $ \ic $ in $ [0,b_1] $.
In $ [b_1,b_2] $, $ f $ is concave. Using the concave version of Theorem~1 (2),
we see that $ f/x $ can have three possible behavior. Since before $ b_1 $,
it is $ \ic $, it will continue to $ \ic $ at least for a little while after $ b_1 $.
This rules out the case (2) i), that it is $ \dc $ in the entire subinterval. 
In the remaining two possibilities,
there may or may not exist a turning point $ c_1\in(b_1,b_2) $, depending
on whether a tangent line through the origin can be drawn touching the 
curve inside $ (b_1,b_2) $. If such a $ c_1 $ exists, then
$ f/x $ will change monotonicity once in $ [0,b_2] $. Otherwise, $ f/x $
remains $ \ic $ in $ [0,b_2] $.

To summarize, in $ [0,b_2] $, $ f/x $ cannot oscillate more
times than $ f' $. Furthermore, the turning point of $ f/g $, $ c_1 $
(if it exists) lags behind that of the latter, namely, $ b_1 $. In other
words, $ c_1>b_1 $.

We can continue with similar arguments in subsequent intervals and it is 
easy to see that the first statement in the above summary remains true
throughout the entire interval $ [0,b) $. The second statement has to be interpreted
in the following way. In each subinterval $ (b_i,b_{i+1}) $, there is at most
one turning point $ c_i $ of $ f'/x $. The sense of the change of monotonicity
(i.e. from $ \ic $ to $ \dc $, or from $ \dc $ to $ \ic) $ is the same as that of $ f' $
at $ b_i $. In the list of all turning points of $ f(x)/x $, some $ c_i $ can be
missing.

The case when $ f(0)>0 $ is just a little more complicated. In view of
Theorem~1 (2), we have to append the possibility of a
first turning point $ c_0\in(0,b_1) $ at which $ f/x $ switches from $ \dc $
to $ \ic $; $ f/x $ can change monotonicity at most $ n+1 $ times.
If there is no such $ c_0 $, then either $ f/x $ is $ \ic $ in $ [0,b_1) $
and the behavior is exactly the same as the case when $ f(0)\leq 0 $,
or $ f/x $ is $ \dc $ in $ [0,b_2] $ and in the remaining subintervals the behavior
mimics that of the case $ f(0)\leq 0 $. In the first situation, $ f/x $ changes
monotonicity at most $ n $ times. In the second situation, it changes at
most $ n-1 $ times.

The case with $ g(x)=x+\gamma  $, $ \gamma >0 $ can be analyzed in a similar way.
After translating to the more general setting, we derive the following
generalization of Pinelis' oscillation rules.

\begin{THEOREM}
Suppose that $ g>0 $ and $ \ic $ in $ [0,b\,] $.
Suppose that $ \phi (u)=f(g^{-1}(u)) $, as defined before,
is a continuous function of $ u $ and is alternatively 
strictly convex $($concave$)$ and strictly concave $($convex$)$ 
in the $ (n+1) $ subintervals corresponding to the decomposition of 
$ [0,b) $ described above. 

If $ g(0)=0 $ and $ f(0)\leq (\geq )\,0 $, then $ f/g $ is initially
$ \ic $ $(\dc)$ in $ [0,b_1) $, while 
in each subsequent subinterval $ (b_i,b_{i+1}) $, $ f/g $ can 
change monotonicity at most once $($in the same sense as the change
of monotonicity of $ f'/g' $ at $b_i )$. Hence, $ f/g $ can oscillate at most $ n $
times.

If $ g(0)=0 $ and $ f(0)>(<)\,0 $, or if $ g(0)>0 $,
then $ f/g $ may or may not have one additional change of 
monotonicity in $ [0,b_1) $. The behavior in subsequent subintervals is the same
as in the previous case. Hence,
$ f/g $ can change monotonicity at most $ n+1 $ times.
\end{THEOREM}

\begin{REM} \rm
Theorem~3 has an obvious analog for + $ \dc $ $ g $.
\end{REM}

Two corollaries have found applications (to be described in the next
section) in some recent work of the author.

\begin{COR}
Suppose $ f $ and $ g $ are continuous in $ [a,b] $, differentiable, with $ g,g'>0 $ in $ (a,b) $.
If $ f'/g' $ is initially $ \ic(\dc) $ and changes monotonicity only once 
in $ [a,b] $, then $ f/g $ has a unique global maximum $($minimum$)$
in $ [a,b] $.
\end{COR}

\begin{COR}
Let $ f $ and $ g $ be continuous and differentiable with $ g,g'>0 $ in $ (a,b) $.
Suppose that $ f'/g' $ is initially $ \ic $ and changes monotonicity once
in $ [a,b] $.
If in addition $ f'(a)g(a)\geq f(a)g'(a) $
and $ f'(b)g(b)\leq f(b)g'(b) $, then $ f/g $ is $ \ic $ in $ [a,b] $
\end{COR}  \begin{PROOF}
Suppose that $ f'/g' $ is $ \ic $ in $ [a,b_1] $ and $ \dc $ in $ [b_1,b] $.
By Corollary 1 (2) i),
the boundary condition at $ x=a $ implies that $ f/g $ is $ \ic $ in $ [a,b_1] $.
In $ [b_1,b] $, we use the concave version of Corollary 1 (2) ii),
applied to the boundary condition at $ b $, to
conclude that $ f/g $ is also $ \ic $ there.
\end{PROOF}

We can push Corollary 4 a little further to still get $ f/g $ $ \ic $
when $ f'/g' $ changes monotonicity two times.

\begin{COR}
Let $ f $ and $ g $ be continuous and differentiable with $ g,g'>0 $ in $ (a,b) $.
Suppose that $ f'/g' $ is initially $ \ic $ and changes monotonicity exactly
twice in $ [a,b] $. Let $ b_2\in(a,b) $ be the second turning point of $ f'/g' $.
If $ f'(a)g(a)\geq f(a)g'(a) $
and $ f'(b_2)g(b_2)\leq f(b_2)g'(b_2) $, then $ f/g $ is increasing in $ [a,b] $
\end{COR}

\begin{REM} \rm
Corollaries 4 and 5 have analogs for proving $ f/g\dc $.
\end{REM}

\section{Examples}

\noindent
\underline{$\vphantom{y}$Example 1}. In some recent work with H. Alzer (on studying some
properties of the error function) we need to compute the global maximum
value of the function
\begin{equation}  k_1(x) = \frac{h(x^2)}{h(x)}  \end{equation}
in $ [0,\infty ) $, where
\begin{equation}  h(x) = \int_{0}^{x} \me^{-t^2}\,dt  \end{equation}
is a multiple of the error function.

Any numerical software can easily produce the estimate
$1.0541564714695\cdots$, attained at $x = 1.246574335142\cdots$.

From a theoretical viewpoint, no matter how accurate the maximization
algorithm is, these values cannot be simply taken to be the correct ones.
That is
because most algorithms can only guarantee to return a local maximum, which
is not necessarily the global maximum. Before any further 
justification, the best we can conclude is that the computed value 
represents a local maximum out of possibly multiple local maxima.
Hence, it
can only be taken as a lower bound of the true value sought. To put any
doubt to rest, we need to
affirm that $ k_1(x) $ changes monotonicity only once.

A first attempt is to show that $ k_1'(x) $ changes sign only
once. Plotting its graph seems to support the claim. Yet no easy
proof is apparent.

Letting $ f(x)=h(x^2) $, we compute the ``H\^opital derivative'' (for lack of
a better name) of $ k_1(x) $
\begin{equation}  \xi (x) = \frac{f'(x)}{h'(x)} = 2x\,\me^{x^2-x^4}  \end{equation}
the derivative of which is
\begin{equation}  \xi '(x) = 2\,\me^{x^2-x^4}(1+2x^2-4x^4) .  \end{equation}
It is easy to verify that $ \xi '(x) $ has a unique positive root $ b_1=\sqrt{\sqrt5+1}/2 $,
and that $ \xi (x) $ is $ \ic $ in $ (0,b_1) $ and $ \dc $ in $ (b_1,\infty ) $.
By Corollary 2, we conclude that $ k_1=f/h $ also
changes monotonicity only once in the same sense, just as desired.

\par\vspace*{1.5\baselineskip}\par
\noindent
\underline{$\vphantom{y}$Example 2}. We also need to know the global maximum of
\begin{equation}  k_2(x) =  \frac{k_1(x)}{x} = \frac{h(x^2)}{xh(x)}  \end{equation}
in $ [0,\infty ) $. The numerical estimate is
$1.0785966957414\cdots$, attained at $x = .68355125808421\cdots$
and we need to ensure that $ k_2(x) $ has only one local maximum.

The case of $ k_2(x) $ is a bit more complicated. With $ f(x)=h(x^2) $
and $ g(x)=xh(x) $, its H\^opital derivative is
\begin{equation}  \xi _1(x) =  \frac{f'(x)}{g'(x)} = \frac{4x\,\me^{-x^4}}{2x\,\me^{-x^2}+h(x)} \,.  \end{equation}
It suffices to show that $ \xi _1(x) $ is initially increasing and then changes
monotonicity only once in $ (0,\infty ) $.

One is tempted to apply the oscillation
rule one more time by computing
\begin{equation}  \xi _2(x) =  \frac{(4x\,\me^{-x^4})'}{(2x\,\me^{-x^2}+h(x))'} = \frac{(4x^4-1)\me^{x^2-x^4}}{x^2-1}  \,.  \end{equation}
However, the conditions of the rule are not satisfied because the denominator of
$ \xi _1(x) $ is not a monotone function of $ x $. In fact, it increases in $ (0,1) $
and decreases in $ (1,\infty ) $. We have to investigate the behaviors of $ \xi _1(x) $
in these two subintervals separately.

First we study the monotonicity of $ \xi _2(x) $ in $ [0,1) $ and $ (1,\infty ) $.
The derivative of $ \xi _2(x) $, after some simplification, is
\begin{equation}  \xi _2'(x) = \frac{2x\,\me^{x^2-x^4}}{(x^2-1)^2} \, (2-11x^2+2x^4+12 x^6-8x^8) \,.  \end{equation}
Since the fraction on the righthand side is nonnegative, the monotonicity of
$ \xi _2(x) $ depends on the sign of the polynomial in parentheses. Letting $ x^2=y $,
the polynomial can be written as
\begin{equation}  p(y) = 2-11y+2y^2+12y^3-4y^4.  \end{equation}
Since $ p(0)=2 $ and $ p(1)=-3 $, $ p(y) $ has a root $ \sigma  $ in $ (0,1) $. We claim
that this is the only positive root, by showing that $ p(y) $ is strictly decreasing
for $ y>0 $. To this end, we note that
$ p'(y)=-11+4y+36y^2-32y^3 $ attains its global maximum in $ (0,\infty ) $ when
$ p''(y)=4+72y-96y^2=0 $. It is easy to verify that this global maximum 
is negative.

It follows that $ \xi _2(x) $ is $ \ic $ in $ (0,\sqrt\sigma ) $ and $ \dc $ in $ (\sqrt\sigma ,1)\cup(1,\infty ) $.

We are now ready to study $ \xi _1(x) $ in $ [1,\infty ) $. Since its denominator
is decreasing in $ [1,\infty ) $, we invoke Corollary 2 instead of Corollary~1.
Since $ \xi _2(x) $ is decreasing, the function $ \phi (u) $ in the hypotheses is concave and we
have to use the concave version of Corollary~2. It is easy to verify that
\begin{equation}  \lim_{x\rightarrow \infty }  \frac{\xi _1(x)}{\xi _2(x)} = \lim_{x\rightarrow \infty } \frac{4x(x^2-1)\me^{-x^2}}{(2x\,\me^{-x^2}+h(x))(4x^4-1)} = 0 .  \end{equation}
Hence, for $ b $ very large, $ \xi _1(b)<\xi _2(b) $. By Corollary~2 (2) ii),
we conclude that $ \xi _1\dc $ in $ [1,b] $ for large $ b $.

In $ [0,1] $, $ \xi _2 $ changes monotonicity once, implying that $ \xi _1 $ changes
monotonicity at most once. The only way that this is compatible with
$ \xi _1\dc $ in $ (1,\infty ) $ is that $ \xi _1 $ changes monotonicity exactly once in
$ [0,\infty ) $, as desired.

\par\vspace*{1.5\baselineskip}\par
\noindent
\underline{$\vphantom{y}$Example 3}. The function
\begin{equation}  k_3(x) = \frac{h(x)-x\,\me^{-x^2}}{x^2}  \end{equation}
occurs in the same study. We want to show that it is $ \ic $ in the interval
$ I=[0,0.967857163] $.
Note that we cannot extend the claim to $ [0,1] $ because $ k_3 $ is
not $ \ic $ at $ x=1 $. Letting $ f(x)=h(x)-x\,\me^{-x^2} $
and $ g(x)=x^2 $, we find that
\begin{equation}  \frac{f'(x)}{g'(x)} = x\,\me^{-x^2} \,.  \end{equation}
However, $ f'/g' $ is not monotone in $ I $;
it changes monotonicity at
$ c=1/\sqrt2 $. Hence, the regular monotone rule fails. We can easily verify
that the hypotheses of Corollary 4 are satisfied and thus conclude that $ k_3 $ is
$ \ic $ in the interval.

Figure 6 shows the graphs of $ k_3(x) $ and its H\^opital derivative. It should
be compared with Figures 3 and 4. The green
dashed curve is not monotone, but as long as it stays above the red
curve, the latter is $ \ic $.

\begin{center}
\FG{50}{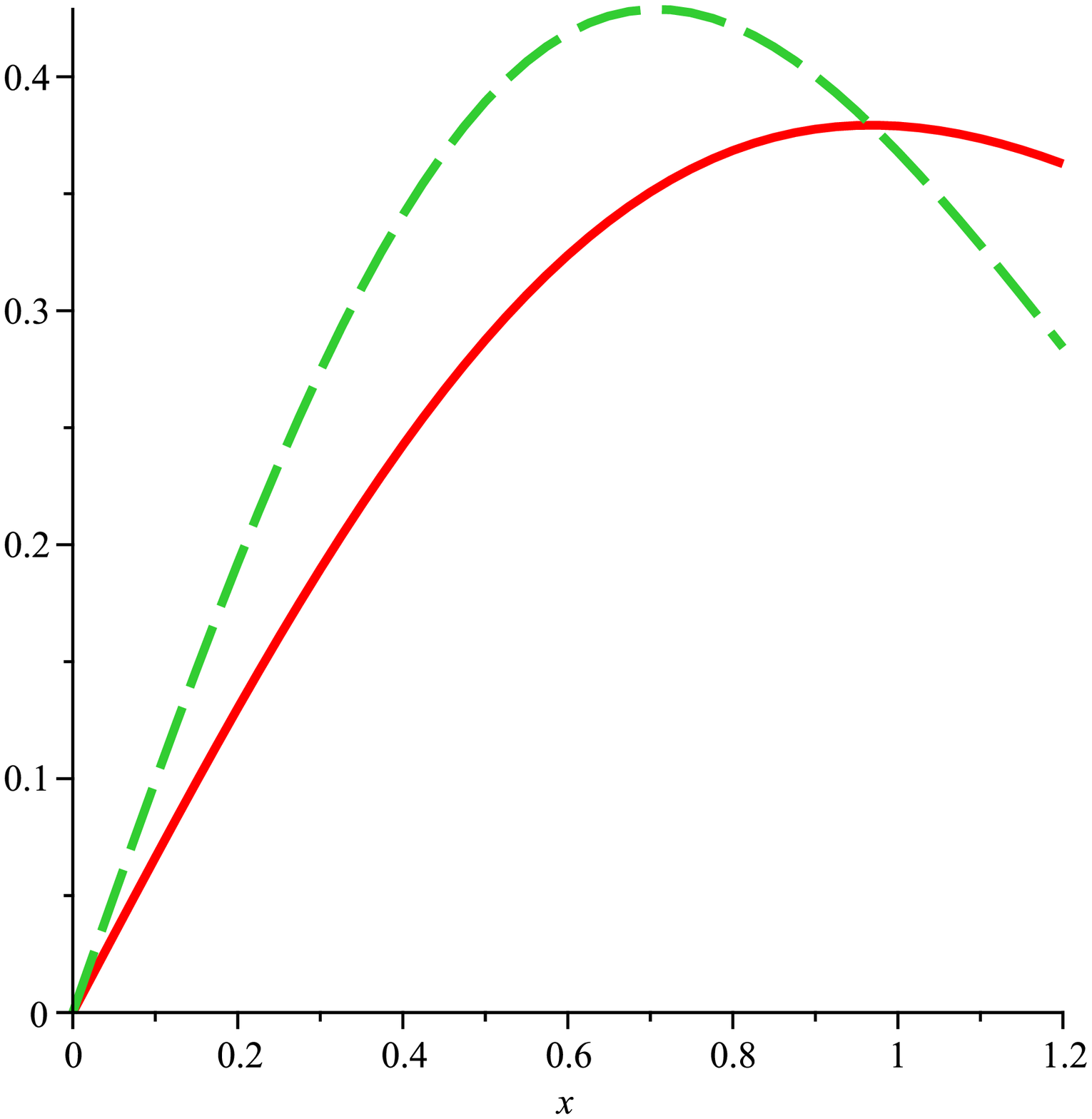}

\par\vspace*{2mm}\par
Figure 6. The graphs of $ k_3(x) $ and its H\^opital derivative.
\end{center}

\par\vspace*{1.5\baselineskip}\par
\noindent
\underline{$\vphantom{y}$Example 4}. The function
\begin{equation}  k_4(x) = \frac{(2x^2-1)h(x)}{h(x)-x\,\me^{-x^2}}   \end{equation}
is $ \ic $ in $ [0,\infty ) $.

The H\^opital derivative is
$$  \xi _3(x) = \frac{(2x^2-1)\,\me^{-x^2}+2xh(x)}{2x^2\me^{-x^2}} \,.  $$
Note that at $ x=0 $, the denominator becomes 0 while the numerator is $ -1 $.
Therefore, the usual monotone rule cannot be used. Instead of using the 
extended rule established in this article, an easier way is to note that
$$  \xi _3(x) = 1 - \frac{1}{2x^2} + \frac{h(x)}{x\,\me^{-x^2}} \,.  $$
The first two terms combined is $ \ic $. Thus, it suffices to show that the
last term is $ \ic $. Yet the usual rule still cannot be used directly
because the denominator of the last term is not monotone in $ [0,\infty ) $.
We can overcome that obstacle by showing that its reciprocal
$ x\,\me^{-x^2}/h(x) $ is $ \dc $. Then the usual rule can be applied.

\end{document}